\documentclass[11pt]{amsart}
\usepackage{latexsym}
\usepackage{amssymb}
\usepackage{amsthm}
\usepackage{tikz}
\usepackage{amscd}

\makeatletter
\@namedef{subjclassname@2010}{\textup{2010} Mathematics Subject Classification}
\makeatother
\newtheorem{theorem}{Theorem}[section]

\newtheorem{proposition}[theorem]{Proposition}

\theoremstyle{remark}
\newtheorem{remark}[theorem]{Remark}
\newtheorem*{acknowledgement}{Acknowledgement}
\theoremstyle{maintheorem}
\newtheorem*{maintheorem}{Main Theorem}

\begin{document}

\title[Orbit spaces of free involutions on lens spaces]{Cohomology algebra of orbit spaces of free involutions on lens spaces}
\author{Mahender Singh}
\address{Indian Institute of Science Education and Research (IISER) Mohali, Phase 9, Sector 81, Knowledge City, S A S Nagar (Mohali), Post Office Manauli, Punjab 140306, India}
\email{mahender@iisermohali.ac.in, mahen51@gmail.com}

\begin{abstract}
Let $G$ be a group acting continuously on a space $X$ and let $X/G$ be its orbit space. Determining the topological or cohomological type of the orbit space $X/G$ is a classical problem in the theory of transformation groups. In this paper, we consider this problem for cohomology lens spaces. Let $X$ be a finitistic space having the mod 2 cohomology algebra of the lens space $L_p^{2m-1}(q_1,\dots,q_m)$. Then we classify completely the possible mod 2 cohomology algebra of orbit spaces of arbitrary free involutions on $X$. We also give examples of spaces realizing the possible cohomology algebras. In the end, we give an application of our results to non-existence of $\mathbb{Z}_2$-equivariant map $\mathbb{S}^n \to X$.
\end{abstract}
\maketitle

\section{Introduction}\label{section1}
Lens spaces are odd dimensional spherical space forms described as follows. Let $p \geq 2$ be a positive integer and $q_1, q_2,\dots, q_m$ be integers coprime to $p$, where $m \geq 1$. Let $\mathbb{S}^{2m-1} \subset {\mathbb{C} }^m$ be the unit sphere and let $\iota^2= -1$. Then the map $$(z_1,\dots,z_m) \mapsto (e^{\frac{2 \pi \iota q_1}{p}}z_1,\dots, e^{\frac{2 \pi \iota q_m}{p}}z_m)$$ defines a free action of the cyclic group $\mathbb{Z}_p$ on $\mathbb{S}^{2m-1}$. The orbit space is called a lens space and is denoted by $L_p^{2m-1}(q_1,\dots,q_m)$. It is a compact Hausdorff orientable manifold of dimension $(2m-1)$. Lens spaces are objects of fundamental importance in topology, particularly in low dimensional topology and have been very well studied. The purpose of this paper is to study these spaces in the context of topological transformation groups.
 
Let $G$ be a group acting continuously on a space $X$ with orbit space $X/G$. Determining the orbit space $X/G$ is a classical problem in the theory of transformation groups. Determining the homeomorphism or homotopy type of the orbit space is often difficult and hence a weaker problem is to determine the possible cohomology algebra of the orbit space.  Orbit spaces of free actions of finite groups on spheres have been studied extensively by Livesay \cite{Livesay}, Rice \cite{Rice}, Ritter \cite{Ritter}, Rubinstein \cite{Rubinstein} and many others. However, not much is known if the space is a compact manifold other than a sphere. Tao \cite{Tao} determined orbit spaces of free involutions on $\mathbb{S}^1 \times \mathbb{S}^2$. Later Ritter \cite{Ritter2} extended the results to free actions of cyclic groups of order $2^n$. Dotzel and others \cite{Dotzel} determined the cohomology algebra of orbit spaces of free $\mathbb{Z}_p$ ($p$ prime) and $\mathbb{S}^1$ actions on cohomology product of two spheres. Recently, the author \cite{msingh1} determined the cohomology algebra of orbit spaces of arbitrary free involutions on cohomology product of two real or complex projective spaces. This was extended to product of finitely many projective spaces by Ashraf \cite{Ashraf}. In this paper, we consider this problem for cohomology lens spaces. Throughout the paper the group actions will be assumed to be continuous.

The mod $p$ cohomology algebra of orbit spaces of free $\mathbb{Z}_p$-actions on cohomology lens spaces was determined recently in \cite{tbsingh} for odd primes $p$. We consider the case of free involutions on cohomology lens spaces. Orbit spaces of involutions on three dimensional lens spaces have been investigated in the literature \cite{Kim1, Myers}. For the lens space $L_p^{3}(q)$, where $p= 4k$ for some $k$, Kim \cite{Kim1} showed that the orbit space of any sense-preserving free involution on $L_p^{3}(q)$ is the lens space $L_{2p}^{3}(q')$, where $ q'q \equiv \pm 1$ or $q' \equiv \pm q$ mod $p$. Myers \cite{Myers} showed that every free involution on a three dimensional lens space is conjugate to an orthogonal free involution, in which case the orbit space is again a lens space (see remark \ref{remark4}).

Let $X \simeq_2 L_p^{2m-1}(q_1,\dots,q_m)$ mean that there is an abstract isomorphism of graded algebras $$H^*(X; \mathbb{Z}_2) \cong H^*(L_p^{2m-1}(q_1,\dots,q_m); \mathbb{Z}_2).$$ We call such a space a mod 2 cohomology lens space and refer to dimension of $L_p^{2m-1}(q_1,\dots,q_m)$ as its dimension. This class of spaces is a big generalization of the class of spaces homotopy equivalent to lens spaces. Two non-trivial examples of this case are $\mathbb{S}^1 \times \mathbb{C}P^{m-1}$ and the Dold manifold $P(1,m-1)$ corresponding to the case $4 \mid p$ (see section \ref{section6}).

Recall that, a finitistic space is a paracompact Hausdorff space whose every open covering has a finite dimensional open refinement, where the dimension of a covering is one less than the maximum number of members of the covering which intersect non-trivially (the notion of a finitistic spaces is due to Swan \cite{Swan}). It is a large class of spaces including all compact Hausdorff spaces and all paracompact spaces of finite covering dimension. Finitistic spaces are the most suitable spaces for studying cohomology theory of transformation groups on general spaces, because of their compatibility with the \v{C}ech cohomology theory (see \cite{Bredon2} for a detailed account of results). Note that the lens space $L_p^{2m-1}(q_1,\dots,q_m)$ is a compact Hausdorff space and hence is finitistic.

From now onwards, for convenience, we write $L_p^{2m-1}(q)$ for $L_p^{2m-1}(q_1,\dots,q_m)$. We consider free involutions on finitistic mod 2 cohomology lens spaces and classify the possible mod 2 cohomology algebra of orbit spaces. More precisely, we prove the following theorem.

\begin{maintheorem}\label{mainthm}
Let $G= \mathbb{Z}_2$ act freely on a finitistic space $X\simeq_2 L_p^{2m-1}(q)$. Then $H^*(X/G; \mathbb{Z}_2) $ is isomorphic to one of the following graded commutative algebras:
\begin{enumerate}
\item $\mathbb{Z}_2[x]/\langle  x^{2m} \rangle$,\\
where  $deg(x)=1$.
\item $\mathbb{Z}_2[x,y]/\langle  x^2,y^m \rangle$,\\
where  $deg(x)=1$ and $deg(y)=2$.
\item $\mathbb{Z}_2[x,y,z]/\langle  x^3, y^2, z^{\frac{m}{2}} \rangle$,\\
where $deg(x)=1$, $deg(y)=1$, $deg(z)=4$ and $m$ is even.
\item $\mathbb{Z}_2[x,y,z]/\langle  x^4, y^2, z^{\frac{m}{2}}, x^2y \rangle,$\\
where $deg(x)=1$, $deg(y)=1$, $deg(z)=4$ and $m$ is even.
\item $\mathbb{Z}_2[x,y,w,z]/ \langle x^5, y^2, w^2, z^{\frac{m}{4}}, x^2y, wy\rangle$,\\
where $deg(x)=1$, $deg(y)=1$, $deg(w)=3$, $deg(z)=8$ and $4\mid m$.
\end{enumerate}
\end{maintheorem}

Our theorem generalizes the results known for orbit spaces of free involutions on the three dimensional lens space $L_p^3(q)$, to that of the large class of finitistic spaces $X \simeq_2 L_p^{2m-1}(q_1,\dots,q_m)$ (see remarks \ref{remark2}, \ref{remark3} and \ref{remark4}).

The paper is organized as follows. We recall the cohomology of lens spaces in section \ref{section2} and give an example of a free involution in section \ref{section3}. After recalling some preliminary results in section \ref{section4}, we prove the main theorem in a sequence of propositions in section \ref{section5}. In section \ref{section6}, we provide some examples realizing the possible cohomology algebras. Finally, in section \ref{section7}, we give an application to $\mathbb{Z}_2$-equivariant maps from $\mathbb{S}^n \to X$ and also suggest another possible application to parametrized Borsuk-Ulam problem.
\bigskip

\section{Cohomology of lens spaces}\label{section2}
The homology groups of a lens space can be easily computed using its cell decomposition (see for example \cite[p.144]{Hatcher}) and are given by
\begin{displaymath}
H_i(L_p^{2m-1}(q); \mathbb{Z} ) = \left\{ \begin{array}{ll}
\mathbb{Z} & \textrm{if $i= 0, ~ 2m-1$}\\
 \mathbb{Z}_p & \textrm{if $i$ is odd and $0 <i< 2m-1$}\\
 0 & \textrm{otherwise.}
  \end{array} \right.
\end{displaymath}

If $p$ is odd, then the mod 2 cohomology groups are
\begin{displaymath}
H^i(L_p^{2m-1}(q); \mathbb{Z}_2 ) = \left\{ \begin{array}{ll}
\mathbb{Z}_2 & \textrm{if $i= 0,~ 2m-1$}\\
 0 & \textrm{otherwise.}
  \end{array} \right.
\end{displaymath}

And if $p$ is even, then
\begin{displaymath}
H^i(L_p^{2m-1}(q); \mathbb{Z}_2 ) = \left\{ \begin{array}{ll}
\mathbb{Z}_2 & \textrm{if $ 0 \leq i \leq 2m-1$}\\
 0 & \textrm{otherwise.}
  \end{array} \right.
\end{displaymath}
\bigskip

\section{Free involutions on lens spaces}\label{section3}
We now construct a free involution on the lens space $L_p^{2m-1}(q)$. Let $q_1,\dots,q_m$ be odd integers coprime to $p$. Consider the map $\mathbb{C}^m \to \mathbb{C}^m$ given by $$(z_1,\dots,z_m) \mapsto (e^{\frac{2 \pi \iota q_1}{2p}}z_1,\dots, e^{\frac{2 \pi \iota q_m}{2p}}z_m ).$$ This map commutes with the $\mathbb{Z}_p$-action on $\mathbb{S}^{2m-1}$ defining the lens space and hence descends to a map $$\alpha:L_p^{2m-1}(q) \to L_p^{2m-1}(q)$$ such that $\alpha^2$ = identity. Thus  $\alpha$ is an involution. Denote an element of $L_p^{2m-1}(q)$ by $[z]$ for $z = (z_1,\dots,z_m) \in \mathbb{S}^{2m-1}$. If  $\alpha([z])= [z]$, then $$(e^{\frac{2 \pi \iota q_1}{2p}}z_1,\dots, e^{\frac{2 \pi \iota q_m}{2p}}z_m )= (e^{\frac{2 \pi \iota k{q_1}}{p}}z_1,\dots, e^{\frac{2 \pi \iota k{q_m}}{p}}z_m )$$ for some integer $k$. Let $1 \leq i\leq m$ be an integer such that $z_i \neq 0$, then $e^{\frac{2 \pi \iota q_i}{2p}}z_i = e^{\frac{2 \pi \iota k{q_i}}{p}}z_i$ and hence  $e^{\frac{2 \pi \iota q_i}{2p}}= e^{\frac{2 \pi \iota k{q_i}}{p}}$. This implies $$\frac{q_i}{2p}-\frac{k{q_i}}{p}= \frac{q_i(1-2k)}{2p}$$ is an integer, a contradiction. Hence the involution $\alpha$ is free. Observe that the orbit space of the above involution is $L_{p}^{2m-1}(q)/ \langle \alpha \rangle = L_{2p}^{2m-1}(q)$.
\bigskip

\section{Preliminary results}\label{section4}
In this section, we recall some basic facts that we will use in the rest of the paper without mentioning explicitly. We refer to \cite{Bredon2, Mccleary} for the details of most of the content in this section. Throughout we will use \v{C}ech cohomology. Let the group $G=\mathbb{Z}_2$ act on a space $X$. Let $$G \hookrightarrow E_G \longrightarrow B_G$$ be the universal principal $G$-bundle. Consider the diagonal action of $G$ on $X \times E_G$. Let $$X_G=(X \times E_G)/G$$ be the orbit space corresponding to the diagonal action. Then the projection $X \times E_G \to E_G$ is $G$-equivariant and gives a fibration $$X\hookrightarrow X_G \longrightarrow B_G$$ called the Borel fibration \cite[Chapter IV]{Borel3}. We will exploit heavily the Leray spectral sequence associated to the Borel fibration $X \hookrightarrow X_G \longrightarrow B_G$, as given in the following theorem.

\begin{proposition}\cite[Theorem 5.2]{Mccleary}\label{proposition4.1}
Let $X\stackrel{i}{\hookrightarrow} X_G \stackrel{\rho}{\longrightarrow} B_G$ be the Borel fibration associated to a $G$-space $X$. Then there is a first quadrant spectral sequence of algebras $\{{E_r}^{*,*}, d_r \}$, converging to $H^*(X_G; \mathbb{Z}_2)$ as an algebra, with $$ {E_2}^{k,l}= H^k(B_G; \mathcal{H}^l(X; \mathbb{Z}_2)),$$ the cohomology of the base $B_G$ with locally constant coefficients $\mathcal{H}^l(X; \mathbb{Z}_2)$ twisted by a canonical action of $\pi_1(B_G)$.
\end{proposition}

The graded commutative algebra $H^*(X_G; \mathbb{Z}_2)$ is isomorphic to the graded commutative algebra Tot$E_{\infty}^{*,*}$, the total complex of $E_{\infty}^{*,*}$, given by $$(\textrm{Tot}E_{\infty}^{*,*})^q= \bigoplus_{k+l=q}E_{\infty}^{k,l}.$$ If the fundamental group $\pi_1(B_G)= \mathbb{Z}_2$ acts trivially on the cohomology $H^*(X; \mathbb{Z}_2)$, the system of local coefficients is constant by \cite[Proposition 5.5]{Mccleary} and hence we have $${E_2}^{k,l} \cong H^k(B_G;\mathbb{Z}_2) \otimes H^l(X; \mathbb{Z}_2).$$
Also, note that $H^*(X_G)$ is a $H^*(B_G)$-module with the multiplication given by $$(b,x) \mapsto \rho^*(b) x$$ for $b \in H^*(B_G)$ and $x \in H^*(X_G)$. Here the product on the right hand side is the cup product.

\begin{proposition}\cite[Theorem 5.9]{Mccleary}\label{proposition4.2}
The edge homomorphisms
$$H^k(B_G;\mathbb{Z}_2)=E_2^{k,0} \longrightarrow E_3^{k,0}\longrightarrow \cdots  \longrightarrow E_k^{k,0} \longrightarrow E_{k+1}^{k,0}=E_{\infty}^{k,0}\subset H^k(X_G;\mathbb{Z}_2)$$
and $$H^l(X_G;\mathbb{Z}_2) \longrightarrow E_{\infty}^{0,l}= E_{l+1}^{0,l} \subset E_{l}^{0,l} \subset \cdots \subset E_2^{0,l}= H^l(X;\mathbb{Z}_2)$$
are the homomorphisms $$\rho^*: H^k(B_G;\mathbb{Z}_2) \to H^k(X_G;\mathbb{Z}_2) ~ ~ ~ \textrm{and} ~ ~ ~ i^*: H^l(X_G;\mathbb{Z}_2)  \to H^l(X;\mathbb{Z}_2).$$
\end{proposition}

Now we recall some results regarding $\mathbb{Z}_2$-actions on finitistic spaces.

\begin{proposition}\cite[Chapter VII, Theorem 1.5]{Bredon2}\label{proposition4.3}
Let $G = \mathbb{Z}_2$ act freely on a finitistic space $X$. Suppose that $H^j(X; \mathbb{Z}_2) = 0$ for all $j > n$. Then $H^j(X_G ; \mathbb{Z}_2)=0$ for $j > n$.
\end{proposition}

Let $h:X_G \to X/G$ be the map induced by the $G$-equivariant projection $X \times E_G \to X$. Then the following is true.

\begin{proposition}\cite[Chapter VII, Proposition 1.1]{Bredon2}\label{proposition4.4}
Let $G=\mathbb{Z}_2$ act freely on a finitistic space $X$. Then $$h^*: H^*(X/G; \mathbb{Z}_2) \stackrel{\cong}{\longrightarrow} H^*(X_G; \mathbb{Z}_2).$$
\end{proposition}

In fact $X/G$ and $X_G$ have the same homotopy type.

\begin{proposition}\cite[Chapter VII, Theorem 1.6]{Bredon2}\label{proposition4.5}
Let $G=\mathbb{Z}_2$ act freely on a finitistic space $X$. Suppose that $\sum_{i\geq0} rank\big( H^i(X;\mathbb{Z}_2)\big) < \infty$ and the induced action on $H^*(X; \mathbb{Z}_2)$ is trivial. Then the Leray spectral sequence associated to $X \hookrightarrow X_G \longrightarrow B_G$ does not degenerate at the $E_2$ term.
\end{proposition}

Recall that, for $G=\mathbb{Z}_2$, we have $$H^*(B_G; \mathbb{Z}_2) \cong \mathbb{Z}_2[t],$$ where $t$ is a homogeneous element of degree one. From now onwards our cohomology groups will be with $\mathbb{Z}_2$ coefficients and we will suppress it from the cohomology notation.
\bigskip

\section{Proof of the Main Theorem}\label{section5}
This section is divided into three subsections according to the various values of $p$. The main theorem follows from a sequence of propositions proved in this section.

\subsection{When $p$ is odd}
Recall that, for $p$ odd, we have $L_p^{2m-1}(q) \simeq_2 \mathbb{S}^{2m-1}$. It is well known that the orbit space of any free involution on a mod 2 cohomology sphere is a mod 2 cohomology real projective space of same dimension (see for example Bredon \cite[p.144]{Bredon3}). For the sake of completeness, we give a quick proof using the Leray spectral sequence.

\begin{proposition}\label{proposition5.1}
Let $G=\mathbb{Z}_2$ act freely on a finitistic space $X \simeq_2 \mathbb{S}^{n}$, where $n \geq 1$. Then $$H^*(X/G; \mathbb{Z}_2) \cong \mathbb{Z}_2[x]/\langle x^{n+1}\rangle,$$
where $deg(x)=1$.
\end{proposition}

\begin{proof}
Note that $E_{2}^{k,l}$ is non-zero only for $l=0,n$. Therefore the differentials $d_r=0$ for $2 \leq r \leq n$ and for $r \geq n+2$. As there are no fixed points, the spectral sequence does not degenerate and hence $$d_{n+1}:E_{n+1}^{k,n} \to E_{n+1}^{k+n+1,0}$$ is non-zero and it is the only non-zero differential. Thus $E_{\infty}^{*,*} = E_{n+2}^{*,*}$ and
\begin{displaymath}
H^j(X_G)= E_{\infty}^{j,0}=\left\{ \begin{array}{ll}
\mathbb{Z}_2 & \textrm{if $0 \leq j \leq n$ }\\
0 & \textrm{otherwise.}
  \end{array} \right.
\end{displaymath}

Let $x= \rho^*(t) \in E_{\infty}^{1,0} \subset H^1(X_G)$ be determined by $t\otimes 1 \in E_2^{1,0}$. Since the cup product $$x \smile (-) : H^j(X_G) \to H^{j+1}(X_G)$$ is an isomorphism for $0 \leq j \leq n-1$, we have $x^j \neq 0$ for $1 \leq j \leq n$. Therefore $$H^*(X_G) \cong \mathbb{Z}_2[x]/\langle x^{n+1}\rangle,$$ where $deg(x)=1$. As the action of $G$ is free, $H^*(X/G ) \cong H^*(X_G)$. This gives the case (1) of the main theorem.
\end{proof}

\subsection{When $p$ is even and $4 \nmid p$}
Let $p$ be even, say $p=2p'$ for some integer $p' \geq 1$. Since $q_1,\dots,q_m$ are coprime to $p$, all of them are odd. Also all of them are coprime to $p'$. Note that $L_p^{2m-1}(q)= L_{2p'}^{2m-1}(q)= L_{p'}^{2m-1}(q)/\langle \alpha \rangle$, where $\alpha$ is the involution on $L_{p'}^{2m-1}(q)$ as defined in section 3. When $4 \nmid p$, we have $p'$ is odd and hence $L_{p'}^{2m-1}(q)\simeq_2 \mathbb{S}^{2m-1}$. This gives $L_p^{2m-1}(q) \simeq_2 \mathbb{R}P^{2m-1}$. Therefore it amounts to determining the mod 2 cohomology algebra of orbit spaces of free involutions on odd dimensional mod 2 cohomology real projective spaces.

\begin{proposition}\label{proposition5.2}
Let $G=\mathbb{Z}_2$ act freely on a finitistic space $X \simeq_2 \mathbb{R}P^{2m-1}$, where $m \geq 1$. Then $$H^*(X/G; \mathbb{Z}_2) \cong \mathbb{Z}_2[x,y]/\langle x^2, y^m\rangle,$$
where $deg(x)=1$ and $deg(y)=2$.
\end{proposition}

\begin{proof}
Clearly the proposition is obvious for $m=1$. Assume that $m > 1$. Let $a \in H^1(X)$ be the generator of the cohomology algebra $H^*(X)$. As there are no fixed points, the spectral sequence does not degenerate at the $E_2$ term. Therefore $d_2(1\otimes a)=t^2 \otimes 1$. One can see that $$d_2:E_2^{k,l} \to E_2^{k+2,l-1}$$ is the trivial homomorphism for $l$ even and an isomorphism for $l$ odd. This gives
\begin{displaymath}
E_3^{k,l}= \left\{ \begin{array}{ll}
\mathbb{Z}_2 & \textrm{if $ k$ = 0, 1 and $l$ = 0, 2,\dots, $2m-2$}\\
 0 & \textrm{otherwise.}
  \end{array} \right.
\end{displaymath}

Note that $d_r = 0$ for all $r \geq 3$ and for all $k, l$ as $E_r^{k+r, l-r+1}= 0$. Hence $E_{\infty}^{*,*} = E_3^{*,*}$ and
\begin{displaymath}
H^j(X_G)= \left\{ \begin{array}{ll}
E_{\infty}^{0,j} & \textrm{if $j$ even}\\
 E_{\infty}^{1,j-1} & \textrm{if $j$ odd.}
  \end{array} \right.
\end{displaymath}

Therefore the additive structure of $H^*(X_G)$ is given by
\begin{displaymath}H^j(X_G)= \left\{ \begin{array}{ll}
\mathbb{Z}_2 & \textrm{if $0 \leq j \leq 2m-1$}\\
0 & \textrm{otherwise.}
  \end{array} \right.
\end{displaymath}

Let $x= \rho^*(t) \in E_{\infty}^{1,0}$ be determined by $t\otimes 1 \in E_2^{1,0}$ and $x^2 \in E_{\infty}^{2,0}=0$. The element $1\otimes a^2 \in E_2^{0,2}$ is a permanent cocycle and determines an element $y \in E_{\infty}^{0,2}=H^2(X_G)$. Also $i^*(y)= a^2$ and $E_{\infty}^{0,2m}=0$ implies $y^m=0$. Since the cup product by $x$ $$x \smile(-): H^j(X_G ) \to  H^{j+1}(X_G)$$
is an isomorphism for $0 \leq j \leq 2m-2$, we have $xy^j \neq 0$ for $0 \leq j \leq m-1$.
Therefore we have $$H^*(X_G) \cong \mathbb{Z}_2[x,y]/\langle x^2, y^m\rangle,$$ where $deg(x)=1$ and $deg(y)=2$. As the action of $G$ is free, we have $H^*(X/G ) \cong H^*(X_G)$. This is the case (2) of the main theorem.
\end{proof}

\begin{remark}\label{remark1}
It is well known that there is no free involution on a finitistic space $X \simeq_2 \mathbb{R}P^{2m}$. For, by the Floyd's Euler characteristic formula \cite[p.145]{Bredon2} $$\chi(X)+ \chi(X^G)= 2 \chi(X/G),$$ $\chi(X)$ must be even, which is a contradiction.
\end{remark}

\begin{remark}\label{remark2}
The above result follows easily for free involutions on $\mathbb{R}P^3$. Let there be a free involution on $\mathbb{R}P^3$. This lifts to a free action on $\mathbb{S}^3$ by a group $H$ of order 4 and $\mathbb{R}P^3/\mathbb{Z}_2 = \mathbb{S}^3/H$. There are only two groups of order $4$, namely, the cyclic group $\mathbb{Z}_{4}$ and $\mathbb{Z}_2 \oplus\mathbb{Z}_2$. By Milnor \cite{Milnor}, $\mathbb{Z}_2 \oplus\mathbb{Z}_2$ cannot act freely on $\mathbb{S}^3$. Hence $H$ must be the cyclic group $\mathbb{Z}_{4}$. Now by Rice \cite{Rice}, this action is equivalent to an orthogonal free action and hence $\mathbb{R}P^3/ \mathbb{Z}_2 = L_4^{3}(q)$.
\end{remark}

\subsection{When $4\mid p$}\label{section5.3}
As above $L_p^{2m-1}(q)= L_{2p'}^{2m-1}(q)= L_{p'}^{2m-1}(q)/\langle\alpha \rangle$. Since $4 \mid p$, $p'$ is even, the cohomology groups $H^i(L_{p'}^{2m-1}(q)) = \mathbb{Z}_2$ for each $0 \leq i \leq 2m-1$ and 0 otherwise. The Smith-Gysin sequence of the orbit map $\eta:L_{p'}^{2m-1}(q)\to L_{2p'}^{2m-1}(q)$, which is a 0-sphere bundle, is given by

{\setlength\arraycolsep{3pt}
\begin{eqnarray}
\lefteqn{ 0 \to H^0(L_{2p'}^{2m-1}(q))  \stackrel{\eta^*}{\to} H^0(L_{p'}^{2m-1}(q)) \stackrel{\tau}{\to} H^0(L_{2p'}^{2m-1}(q)) \stackrel{\smile u}{\to} H^{1}(L_{2p'}^{2m-1}(q))  \stackrel{\eta^*}{\to} \cdots}\nonumber\\
& & \cdots \stackrel{\smile u}{\to} H^{2m-1}(L_{2p'}^{2m-1}(q)) \stackrel{\eta^*}{\to} H^{2m-1}(L_{p'}^{2m-1}(q)) \stackrel{\tau}{\to} H^{2m-1}(L_{2p'}^{2m-1}(q)) \to 0, \nonumber
\end{eqnarray}}
where $\tau$ is the transfer map. The cup square $u^2$ of the characteristic class $u \in H^1( L_{2p'}^{2m-1}(q))$ is zero by the exactness of the sequence. This gives the cohomology algebra $$H^*(L_{2p'}^{2m-1}(q)) \cong \wedge[u] \otimes \mathbb{Z}_2[v]/\langle v^m\rangle \cong \mathbb{Z}_2[u,v]/\langle u^2, v^m\rangle,$$ where $u \in H^1(L_{2p'}^{2m-1}(q))$ and $v \in H^2(L_{2p'}^{2m-1}(q))$.

Two non-trivial examples of this case are $\mathbb{S}^1 \times \mathbb{C}P^{m-1}$ and the Dold manifold $P(1,m-1)$. We will elaborate these examples in section \ref{section6}.

Let $u, v$ be generators of $H^*(X)= H^*(L_{2p'}^{2m-1}(q))$ as above. As the group $G= \mathbb{Z}_2$ acts freely on $X$ with trivial action on $H^*(X)$, the spectral sequence does not degenerate at the $E_2$ term. If $d_2 = 0$, then $d_3\neq 0$, otherwise, the spectral sequence degenerate at the $E_2$ term. Thus we have the following proposition.

\begin{proposition}\label{proposition5.3}
Let $G= \mathbb{Z}_2$ act freely on a finitistic space $X\simeq_2 L_p^{2m-1}(q)$, where $4 \mid p$. Let $\{E_r^{*,*}, d_r \}$ be the Leray spectral sequence associated to the fibration $X\stackrel{i}{\hookrightarrow} X_G \stackrel{\rho}{\longrightarrow} B_G$. If $u$, $v$ are generators of $H^*(X)$ such that $d_2 = 0$, then $$H^*(X/G; \mathbb{Z}_2) \cong \mathbb{Z}_2[x,y,z]/\langle  x^3, y^2, z^{\frac{m}{2}} \rangle,$$ where $deg(x)=1$, $deg(y)=1$, $deg(z)=4$ and $m$ is even.
\end{proposition}

\begin{proof}
As $d_2 = 0$, we have that $d_3\neq 0$, otherwise, the spectral sequence degenerate at the $E_2$ term. Since $d_3(1\otimes u)=0$, we must have  $d_3(1\otimes v) =t^3 \otimes 1$. By the multiplicative property of $d_3$, we have
\begin{displaymath}
d_3(1 \otimes v^q) = \left\{ \begin{array}{ll}
t^3 \otimes v^{q-1} & \textrm{if $0 < q < m$ odd}\\
0 & \textrm{if $0 < q < m$ even}\\
\end{array} \right.
\end{displaymath}

Similarly
\begin{displaymath}
d_3(1 \otimes uv^q) = \left\{ \begin{array}{ll}
t^3 \otimes uv^{q-1} & \textrm{if $0 < q < m$ odd}\\
0 & \textrm{if $0 < q < m$ even}\\
\end{array} \right.
\end{displaymath}

This shows that $$d_3:E_3^{k,l} \to E_3^{k+3,l-2}$$ is an isomorphism for $l= 4q+2, 4q+3$ and zero for $l= 4q, 4q+1$. Also note that $v^m = 0$. If $m$ is odd, then $$0 = d_3(1 \otimes v^m)= d_3 \big((1 \otimes v^{m-1})(1 \otimes v)\big)= t^3 \otimes v^{m-1},$$ a contradiction. Hence $m$ must be even, say $m=2n$ for some $n \geq 1$. Therefore
\begin{displaymath}
E_4^{k,l} = \left\{ \begin{array}{ll}
E_3^{k,l} & \textrm{if $k=0,1,2$ and $l=4q, 4q+1$, where  $0 \leq q \leq n-1$}\\
0 & \textrm{otherwise.}
  \end{array} \right.
\end{displaymath}

Note that $d_r=0$ for all $r \geq 4$ and for all $k,l$ as $E_r^{k+r,l-r+1}= 0$. Therefore $E_{\infty}^{*,*}=E_4^{*,*}$ and the additive structure of $H^*(X_G)$ is given by
\begin{displaymath}H^j(X_G)= \left\{ \begin{array}{ll}
\mathbb{Z}_2 & \textrm{if $j=4q, 4q+3$, where  $0 \leq q \leq n-1$}\\
\mathbb{Z}_2 \oplus \mathbb{Z}_2 & \textrm{if $j=4q+1, 4q+2$, where  $0 \leq q \leq n-1$}\\
0 & \textrm{otherwise.}
\end{array} \right.
\end{displaymath}

Let $x= \rho^*(t) \in E_{\infty}^{1,0}$ be determined by $t \otimes 1 \in E_{2}^{1,0}$. As $E_{\infty}^{3,0}=0$, we have $x^3=0$. Note that $1 \otimes u \in E_{2}^{0,1}$ is a permanent cocycle and hence determines an element say $y \in E_{\infty}^{0,1}$. Also $i^*(y)=u$ and $E_{\infty}^{0,2}=0$ implies $y^2=0$. Similarly $1 \otimes v^2$ is a permanent cocycle and therefore it determines an element say $z \in E_{\infty}^{0,4}= H^4(X_G)$. Also $i^*(z)=v^2$ and $E_{\infty}^{0,4n}=0$ implies $z^{\frac{m}{2}}=0$. Since the cup product by $x$ $$x \smile(-): H^j(X_G ) \to  H^{j+1}(X_G)$$ is an isomorphism for $0 \leq j \leq 2m-2$, we have $xz^r \neq 0$ for $0 \leq r \leq \frac{m-1}{2}$. Therefore $$H^*(X_G) \cong \mathbb{Z}_2[x,y,z]/\langle  x^3, y^2, z^{\frac{m}{2}} \rangle,$$ where $deg(x)=1$, $deg(y)=1$ and $deg(z)=4$. As the action of $G$ is free, we have $H^*(X/G ) \cong H^*(X_G)$. This is the case (3) of the main theorem.
\end{proof}

Next, we consider $d_2 \neq 0$, for which we have the following possibilities:
\begin{itemize}
\item[(A)] $d_2(1\otimes u)=t^2 \otimes 1$ and $d_2(1\otimes v)=t^2 \otimes u$,
\item[(B)] $d_2(1\otimes u)=t^2 \otimes 1$ and $d_2(1\otimes v)=0$ and
\item[(C)]$d_2(1\otimes u)=0$ and $d_2(1\otimes v)= t^2 \otimes u$.
\end{itemize}

We consider the above possibilities one by one. We first observe that the possibility (A) does not arise. Suppose that $d_2(1\otimes u) =t^2 \otimes 1$ and $d_2(1\otimes v)= t^2 \otimes u$. By the multiplicative property of $d_2$, we have
\begin{displaymath}
d_2(1 \otimes v^q) = \left\{ \begin{array}{ll}
t^2 \otimes uv^{q-1} & \textrm{if $0 < q < m$ odd}\\
0 & \textrm{if $0 < q < m$ even}\\
\end{array} \right.
\end{displaymath}
and $d_2(1 \otimes uv^q)= t^2 \otimes v^q$ for $0 < q < m$. This shows that $$d_2:E_2^{k,l} \to E_2^{k+2,l-1}$$ is an isomorphism if $l$ even and $4 \nmid l$ or $l$ odd. And $d_2$ is zero if $4 \mid l$. Just as in the previous proposition, $m$ must be even, say $m=2n$ for some $n \geq 1$. This gives
\begin{displaymath}
E_3^{k,l} = \left\{ \begin{array}{ll}
E_2^{k,l} & \textrm{if $k=0,1$ and  $l=4q$, where  $0 \leq q \leq n-1$}\\
0 & \textrm{otherwise.}
  \end{array} \right.
\end{displaymath}

Note that $d_r=0$ for all $r \geq 3$ and for all $k,l$ as $E_r^{k+r, l-r+1}= 0$. Therefore $E_{\infty}^{*,*}=E_3^{*,*}$ and
\begin{displaymath}H^j(X_G)= \left\{ \begin{array}{ll}
\mathbb{Z}_2 & \textrm{if $l=4q, 4q+1$, where  $0 \leq q \leq n-1$}\\
0 & \textrm{otherwise.}
  \end{array} \right.
\end{displaymath}

In particular, this shows that $H^{2m-1}(X/G) \cong H^{2m-1}(X_G)=0$. But the Smith-Gysin sequence $$ \cdots \to H^{2m-1}(X/G) \stackrel{\eta^*}{\to} H^{2m-1}(X) \stackrel{\tau}{\to} H^{2m-1}(X/G) \to 0$$ implies that $H^{2m-1}(X)=0$, which is a contradiction. Hence this possibility does not arise.

For the possibility (B), we have the following proposition.

\begin{proposition}\label{proposition5.5}
Let $G= \mathbb{Z}_2$ act freely on a finitistic space $X\simeq_2 L_p^{2m-1}(q)$, where $4 \mid p$. Let $\{E_r^{*,*}, d_r \}$ be the Leray spectral sequence associated to the fibration $X\stackrel{i}{\hookrightarrow} X_G \stackrel{\rho}{\longrightarrow} B_G$. If $u$, $v$ are generators of $H^*(X)$ such that $d_2(1\otimes u) \neq 0$ and $d_2(1\otimes v)= 0$, then $$H^*(X/G; \mathbb{Z}_2) \cong \mathbb{Z}_2[x,y]/\langle x^2,y^m \rangle,$$ where $deg(x)=1$ and $deg(y)=2$.
\end{proposition}

\begin{proof}
Let $d_2(1\otimes u) = t^2 \otimes 1$ and $d_2(1\otimes v)= 0$. Consider $$d_2:E_2^{k,l} \to E_2^{k+2,l-1}.$$ If $l=2q$, then $d_2(t^k \otimes v^q)=0$ and if $l=2q+1$, then $d_2(t^k \otimes uv^q)= t^{k+2} \otimes v^q$ for $0 \leq q \leq m-1$. This gives
\begin{displaymath}
E_3^{k,l} = \left\{ \begin{array}{ll}
E_2^{k,l} & \textrm{if $k=0,1$ and  $l=0,2,\dots,2m-2$}\\
0 & \textrm{otherwise.}
\end{array} \right.
\end{displaymath}

Note that $$d_r:E_r^{k,l} \to E_r^{k+r,l-r+1}$$ is zero for all $r \geq 3$ and for all $k,l$ as $E_r^{k+r,l-r+1}= 0$. This gives $E_{\infty}^{*,*}= E_3^{*,*}$. But
\begin{displaymath}
H^j(X_G)= \left\{ \begin{array}{ll}
E_{\infty}^{0,j} & \textrm{if $j$ even}\\
E_{\infty}^{1,j-1} & \textrm{if $j$ odd.}
\end{array} \right.
\end{displaymath}

Therefore
\begin{displaymath}H^j(X_G)= \left\{ \begin{array}{ll}
\mathbb{Z}_2 & \textrm{if $0 \leq j \leq 2m-1$}\\
0 & \textrm{otherwise.}
\end{array} \right.
\end{displaymath}

Let $x= \rho^*(t) \in E_{\infty}^{1,0}$ be determined by $t \otimes 1 \in E_{2}^{1,0}$. As $E_{\infty}^{2,0}=0$, we have $x^2=0$. Note that $1 \otimes v$ is a permanent cocycle and therefore it determines an element say $y \in E_{\infty}^{0,2}= H^2(X_G)$. Also $i^*(y)=v$ and $E_{\infty}^{0,2m}=0$ implies $y^m = 0$. Since the cup product by $x$ $$x \smile(-): H^j(X_G ) \to  H^{j+1}(X_G)$$ is an isomorphism for $0 \leq j \leq 2m-2$, we have $xy^r \neq 0$ for $0 \leq r \leq m-1$. Therefore $$H^*(X_G) \cong \mathbb{Z}_2[x,y]/\langle  x^2,y^m \rangle,$$ where $deg(x)=1$ and $deg(y)=2$. As the action of $G$ is free, $H^*(X/G ) \cong H^*(X_G)$. Again we get the case (2) of the main theorem.
\end{proof}

Finally, for the possibility (C), we have the following proposition.

\begin{proposition}\label{proposition5.6}
Let $G= \mathbb{Z}_2$ act freely on a finitistic space $X\simeq_2 L_p^{2m-1}(q)$, where $4 \mid p$. Let $\{E_r^{*,*}, d_r \}$ be the Leray spectral sequence associated to the fibration $X\stackrel{i}{\hookrightarrow} X_G \stackrel{\rho}{\longrightarrow} B_G$. If $u$, $v$ are generators of $H^*(X)$ such that $d_2(1\otimes u) =0$ and $d_2(1\otimes v) \neq 0$, then $H^*(X/G; \mathbb{Z}_2)$ is isomorphic to one of the following graded commutative algbras:
\begin{enumerate}
\item[(i)] $\mathbb{Z}_2[x,y,z]/\langle x^4, y^2, z^{\frac{m}{2}}, x^2y \rangle$,\\
where $deg(x)=1$, $deg(y)=1$, $deg(z)=4$ and $m$ is even.
\item[(ii)] $\mathbb{Z}_2[x,y,w,z]/ \langle x^5, y^2, w^2, z^{\frac{m}{4}}, x^2y, wy\rangle$,\\
where $deg(x)=1$, $deg(y)=1$, $deg(w)=3$, $deg(z)=8$ and $4\mid m$.
\end{enumerate}
\end{proposition}

\begin{proof}
Let $d_2(1\otimes u)=0$ and $d_2(1\otimes v)= t^2 \otimes u$. The derivation property of the differential gives
\begin{displaymath}
d_2(1 \otimes v^q) = \left\{ \begin{array}{ll}
t^2 \otimes uv^{q-1} & \textrm{if $0 < q < m$ odd}\\
0 & \textrm{if $0 < q < m$ even.}\\
\end{array} \right.
\end{displaymath}

Also $d_2(1 \otimes uv^q)= 0$ for $0 < q < m$. Note that $v^m=0$. If $m$ is odd, then $$0 = d_2(1 \otimes v^m)= d_2 \big((1 \otimes v^{m-1})(1 \otimes v)\big)= t^2 \otimes uv^{m-1},$$ a contradiction. Hence $m$ must be even, say $m=2n$ for some $n \geq 1$. From this we get $$d_2:E_2^{k,l} \to E_2^{k+2,l-1}$$ is an isomorphism if $l$ even and $4 \nmid l$, and is zero if $l$ odd or $4 \mid l$. This gives
\begin{equation*}
(\star) \hspace{1cm} E_3^{k,l} = \left\{ \begin{array}{ll}
E_2^{k,l} & \textrm{if $k\geq 0$ arbitrary and $l=4q,4q+3$, where $0\leq q \leq n-1$}\\
E_2^{k,l} & \textrm{if $k=0,1$ and $l=4q+1$, where $0\leq q \leq n-1$}\\
0 & \textrm{otherwise.}
  \end{array} \right.
\end{equation*}

We now consider the differentials one by one. First, we consider $$d_3:E_3^{k,l} \to E_3^{k+3,l-2}.$$ Clearly $d_3=0$ for all $k$ and for $l=4q,4q+3$ as $E_3^{k+3,l-2}=0$ in this case. For $k=0,1$ and for $l=4(q+1)+1= 4q+5$, $$d_3:E_3^{k,4q+5} \to E_3^{k+3,4q+3}$$ is also zero, because if $a \in E_3^{k,4q+5}$ and $d_3(a)=[t^{k+3} \otimes uv^{2q+1}]$, then for $b=[t^2 \otimes 1] \in E_3^{2,0}$, we have $ab \in E_3^{k+2,4q+5}= 0$ and hence $$0 = d_3(ab)= d_3(a)b + ad_3(b)= d_3(a)b + 0 = [t^{k+5} \otimes uv^{2q+1}],$$ which is a contradiction. Hence $d_3=0$ for all $k,l$.

Next we break the remaining proof in the following two cases:
\begin{enumerate}
\item[(a)] $d_4:E_4^{0,3} \to E_4^{4,0}$ is non-zero.
\item[(b)] $d_4:E_4^{0,3} \to E_4^{4,0}$ is zero.
\end{enumerate}

\subsection*{(a) When $d_4:E_4^{0,3} \to E_4^{4,0}$ is non-zero}
Let $d_4([1 \otimes uv])= [t^4 \otimes 1]$. This gives $$d_4:E_4^{k,l} \to E_4^{k+4,l-3}$$ is an isomorphism for all $k$ and for $l= 4q+3$, where $0 \leq q \leq n-1$ and zero otherwise. This gives
\begin{equation*}
E_5^{k,l} = \left\{ \begin{array}{ll}
E_4^{k,l} & \textrm{if $k=0,1,2,3$ and $l=4q$, where $0\leq q \leq n-1$}\\
E_4^{k,l} & \textrm{if $k=0,1$ and $l=4q+1$, where $0\leq q \leq n-1$}\\
0 & \textrm{otherwise.}
  \end{array} \right.
\end{equation*}

It is clear that $d_r=0$ for all $r \geq 5$ and for all $k,l$ as $E_r^{k+r,l-r+1}=0$. Hence, $E_{\infty}^{*,*}= E_5^{*,*}$ and the additive structure of $H^*(X_G)$ is given by
\begin{displaymath}
H^j(X_G)= \left\{ \begin{array}{ll}
\mathbb{Z}_2 & \textrm{if $j = 4q, 4q+3$, where  $0 \leq q \leq n-1$}\\
\mathbb{Z}_2 \oplus \mathbb{Z}_2 & \textrm{if $j= 4q+1, 4q+2$, where  $0 \leq q \leq n-1$}\\
0 & \textrm{otherwise.}
\end{array} \right.
\end{displaymath}

We see that $1 \otimes v^2 \in E_2^{0,4}$ and $1 \otimes u \in E_2^{0,1}$ are permanent cocycles. Hence, they determine elements $z \in E_{\infty}^{0,4}\subseteq H^4(X_G)$ and $y \in E_{\infty}^{0,1}\subseteq H^1(X_G)$, respectively. As $H^4(X_G)= E_{\infty}^{0,4}=E_2^{0,4}$, we have $i^*(z)= v^2$. Since $E_{\infty}^{0,4n}=0$, we get $z^n=0$. Similarly, $i^*(y)=u$ and $E_{\infty}^{0,2}=0$ implies $y^2=0$.

Let $x=\rho^*(t) \in E_{\infty}^{1,0}\subseteq H^1(X_G)$ be determined by $t \otimes 1 \in E_{2}^{1,0}$. As $E_{\infty}^{4,0}=0$, we have $x^4=0$. Also, the cup product  $x^2y \in E_{\infty}^{2,1}=0$. Hence, $$H^*(X_G) \cong \mathbb{Z}_2[x,y,z]/ \langle x^4, y^2, z^{\frac{m}{2}}, x^2y \rangle,$$ where $deg(x)=1$, $deg(y)=1$ and $deg(z)=4$. As the action of $G$ is free, $H^*(X/G ) \cong H^*(X_G)$. This gives the case (4) of the main theorem.

\subsection*{(b) When $d_4:E_4^{0,3} \to E_4^{4,0}$ is zero}
We show that $$d_4:E_4^{k,l} \to E_4^{k+4,l-3}$$ is zero for all $k,l$. Note that $E_4^{k+4,l-3}=0$ for all $k$ and for $l=4q$. Similarly $E_4^{k+4,l-3}=0$ for $k=0,1$ and for $l=4q+1$. Now for any $k$ and $l=4q+3$, we have that $$d_4:E_4^{k,4q+3} \to E_4^{k+4,4q}$$ is given by $$d_4([t^k \otimes uv^{2q+1}])= (d_4[t^k \otimes uv])[1 \otimes v^{2q}]+ [t^k \otimes uv](d_4[1 \otimes v^{2q}])= 0.$$ This shows that $d_4=0$ for all $k,l$.

Now we have the following two subcases:
\begin{enumerate}
\item[(b1)] $d_5:E_5^{0,4} \to E_5^{5,0}$ is non-zero.
\item[(b2)] $d_5:E_5^{0,4} \to E_5^{5,0}$ is zero.
\end{enumerate}

\subsection*{(b1) When $d_5:E_5^{0,4} \to E_5^{5,0}$ is non-zero}
Let $d_5([1 \otimes v^2])= [t^5 \otimes 1]$. Then
\begin{displaymath}
d_5([1 \otimes v^{2q}])=q[t^5 \otimes v^{2(q-1)}] =\left\{ \begin{array}{ll}
[t^5 \otimes v^{2(q-1)}] & \textrm{if $0 < q < m$ odd}\\
0 & \textrm{if $0 < q < m$ even}\\
\end{array} \right.
\end{displaymath}
and
\begin{displaymath}
d_5([1 \otimes uv^{2q+1}]) = q[t^5 \otimes uv^{2q-1}]= \left\{ \begin{array}{ll}
[t^5 \otimes uv^{2q-1}] & \textrm{if $0 < q < m$ odd}\\
0 & \textrm{if $0 < q < m$ even.}\\
\end{array} \right.
\end{displaymath}

Note that $1 \otimes uv^{m+1}=0$ and hence $0=d_5([1 \otimes uv^{m+1}])= n[t^5 \otimes uv^{m-1}].$ But this is possible only when $2 \mid n$ and hence $4 \mid m$.

From above we obtain
\begin{equation*}
E_6^{k,l} = \left\{ \begin{array}{ll}
E_5^{k,l} & \textrm{if $k=0,1,2,3,4$ and $l=8q,8q+3$, where $0\leq q \leq \frac{n-2}{2}$}\\
E_5^{k,l} & \textrm{if $k=0,1$ and $l=8q+1,8q+5$, where $0\leq q \leq \frac{n-2}{2}$}\\
0 & \textrm{otherwise.}
\end{array} \right.
\end{equation*}

One can see that $d_r=0$ for all $r \geq 6$ and for all $k,l$ as $E_r^{k+r,l-r+1}=0$. Hence, $E_{\infty}^{*,*}= E_6^{*,*}$ and the additive structure of $H^*(X_G)$ is given by
\begin{displaymath}
H^j(X_G)= \left\{ \begin{array}{ll}
\mathbb{Z}_2 & \textrm{if $j= 8q,8q+7$, where  $0 \leq q \leq \frac{n-2}{2}$}\\
\mathbb{Z}_2 \oplus \mathbb{Z}_2 & \textrm{if $8q < j < 84q+7$, where  $0 \leq q \leq \frac{n-2}{2}$}\\
0 & \textrm{otherwise.}
\end{array} \right.
\end{displaymath}

We see that $1 \otimes v^4 \in E_2^{0,8}$, $1 \otimes uv \in E_2^{0,3}$ and $1 \otimes u \in E_2^{0,1}$ are permanent cocycles. Hence, they determine elements $z \in E_{\infty}^{0,8}\subseteq H^8(X_G)$, $w \in E_{\infty}^{0,3}\subseteq H^3(X_G)$ and $y \in E_{\infty}^{0,1}\subseteq H^1(X_G)$, respectively. As $H^8(X_G)= E_{\infty}^{0,8}=E_2^{0,8}$, we have $i^*(z)= v^4$. Also $E_{\infty}^{0,2m}=0$ implies $z^{\frac{m}{4}}=0$. Similarly, $i^*(w)=uv$ and $E_{\infty}^{0,6}=0$ implies $w^2=0$. Finally, $i^*(y)=u$ and $E_{\infty}^{0,2}=0$ implies $y^2=0$.

Let $x = \rho^*(t) \in E_{\infty}^{1,0}\subseteq H^1(X_G)$ be determined by $t \otimes 1 \in E_{2}^{1,0}$. As $E_{\infty}^{5,0}=0$, we have $x^5=0$. Also, the only trivial cup products are $x^2y \in E_{\infty}^{2,1}=0$ and $wy \in E_{\infty}^{0,4}=0$. Hence, $$H^*(X_G) \cong \mathbb{Z}_2[x,y,w,z]/ \langle x^5, y^2, w^2, z^{\frac{m}{4}}, x^2y, wy\rangle,$$ where $deg(x)=1$, $deg(y)=1$, $deg(w)=3$ and $deg(z)=8$. As the action of $G$ is free, $H^*(X/G ) \cong H^*(X_G)$. This gives the case (5) of the main theorem.

\subsection*{(b2) When $d_5:E_5^{0,4} \to E_5^{5,0}$ is zero}
We show that $$d_5:E_5^{k,l} \to E_5^{k+5,l-4}$$ is zero for all $k,l$. Note that $E_5^{k+5,l-4}=0$ for $k=0,1$ and for $l=4q+1$. For any $k$ and for $l=4q,4q+3$, $d_5$ is zero by the derivation property of $d_5$ and the above condition (b2). Hence $d_5=0$ for all $k,l$.

We have that $d_2$, $d_3$, $d_4$ and $d_5$ are all zero for the values of $k$ and $l$ given by equation ($\star$). Note that only $1 \otimes u$, $1 \otimes v^2$ and $1 \otimes uv$ survives to the $E_6$ term. For $r \geq 6$, a typical non-zero element in $E_r^{k,l}$ is of the form $[t^k \otimes v^{2q}]$, $[t^k \otimes uv^{2q+1}]$ or $[t^k \otimes uv^{2q}]$ according as $l=4q,4q+3,4q+1$ for $0\leq q \leq n-1$, respectively. But all these elements can be written as a product of previous three elements for which $d_r=0$ for $r \geq 6$. Hence $E_{\infty}^{*,*}= E_3^{*,*}$. This gives  $H^j(X_G) \neq 0$ for $j > 2m-1$ (in particular $H^{2m}(X_G)\neq 0$), which is a contradiction by proposition \ref{proposition4.3}. Hence (b2) does not arise.
\end{proof}

With this we have completed the proof of the main theorem.
\bigskip

\section{Examples realizing the cohomology algebras}\label{section6}
In this section we provide some examples realizing the possible cohomology algebras of the main theorem.

\begin{itemize}
\item The case (1) of the main theorem can be realized by taking any free involution on a sphere.

\item The Smith-Gysin sequence shows that the example discussed in section \ref{section3} realizes the case (2). For another example, let $X=\mathbb{S}^1 \times \mathbb{C}P^{m-1}$, where $m \geq 2$. The mod 2 cohomology algebra of $X$ is given by $$H^*(X) \cong \mathbb{Z}_2[u,v]/\langle u^2, v^m \rangle,$$ where $deg(u)=1$ and $deg(v)=2$. Note that $X$ always admits a free involution as $\mathbb{S}^1$ does so. Taking any free involution on $\mathbb{S}^1$ and the trivial action on $\mathbb{C}P^{m-1}$ gives $X/G=\mathbb{S}^1 \times \mathbb{C}P^{m-1}$. Hence $$H^*(X/G) \cong \mathbb{Z}_2[x,y]/ \langle x^2, y^m \rangle,$$ where $deg(x)=1$ and $deg(y)=1$. This also realizes the case (2) of the main theorem.

\item We now construct an example for the case (3). Let $X$ be as above. If $m$ is even, then $\mathbb{C}P^{m-1}$ always admits a free involution. In fact, if we denote an element of $\mathbb{C}P^{m-1}$ by $[z_1,z_2,\dots,z_{m-1},z_m]$, then 
the map $$[z_1,z_2,\dots,z_{m-1},z_m] \mapsto [-\overline{z}_2,\overline{z}_1,\dots,-\overline{z}_m,\overline{z}_{m-1}]$$ defines an involution on $\mathbb{C}P^{m-1}$. Now if $$[z_1,z_2,\dots,z_{m-1},z_m]=[-\overline{z}_2,\overline{z}_1,\dots,-\overline{z}_m,\overline{z}_{m-1}],$$ then there exits a $\lambda \in \mathbb{S}^1$ such that $$(\lambda z_1,\lambda z_2,\dots,\lambda z_{m-1},\lambda z_m)=(-\overline{z}_2,\overline{z}_1,\dots,-\overline{z}_m,\overline{z}_{m-1}).$$ This gives $z_1=z_2=\dots =z_{m-1}=z_m=0$, a contradiction. Hence the involution is free. The mod 2 cohomology algebra of orbit spaces of free involutions on odd dimensional complex projective spaces was determined by the author in \cite[Corollary 4.7]{msingh1}. More precisely, it was proved that: For any free involution on $\mathbb{C}P^{m-1}$, where $m \geq 2$ is even, the mod 2 cohomology algebra of the orbit space is given by $$H^*(\mathbb{C}P^{m-1}/G) \cong \mathbb{Z}_2[x,z]/ \langle x^3, z^{\frac{m}{2}}\rangle,$$ where $deg(x)=1$ and $deg(z)=4$. Taking the trivial involution on $\mathbb{S}^1$ and a free involution on $\mathbb{C}P^{m-1}$, we have that $X/G=\mathbb{S}^1 \times \big( \mathbb{C}P^{m-1}/G \big)$. Using the above result, we have $$H^*(X/G) \cong \mathbb{Z}_2[x,y,z]/ \langle x^3, y^2, z^{\frac{m}{2}}\rangle,$$ where $deg(x)=1$, $deg(y)=1$ and $deg(z)=4$. This realizes the case (3) of the main theorem.

\item We do not have examples realizing the cases (4) and (5) of the main theorem. However, we feel that Dold manifolds may give some examples realizing these cases. For integers $r,s \geq 0$, a Dold manifold $P(r,s)$ is defined as $$P(r,s)= \mathbb{S}^r \times \mathbb{C}P^s/ \sim,$$ where $\big((x_1,\dots,x_{r+1}),[z_1,\dots,z_{s+1}]\big) \sim \big((-x_1,\dots,-x_{r+1}),[\overline{z}_1,\dots,\overline{z}_{s+1}]\big)$. Consider the equivariant projection $\mathbb{S}^r \times \mathbb{C}P^s \to \mathbb{S}^r$. On passing to orbit spaces, the Dold manifold can also be seen as the total space of the fiber bundle $$\mathbb{C}P^s \hookrightarrow P(r,s) \to \mathbb{R}P^r.$$ The mod 2 cohomology algebra of a Dold manifold is well known \cite{Dold1} and is given by $$H^*(P(r,s); \mathbb{Z}_2) \cong \mathbb{Z}_2[x,y]/\langle x^{r+1}, y^{s+1} \rangle,$$ where $deg(x)=1$ and $deg(y)=2$. Note that the Dold manifold $P(1,m-1)\simeq_2 L_p^{2m-1}(q)$ for $4\mid p$ and can be considered as the twisted analogue of $X=\mathbb{S}^1 \times \mathbb{C}P^{m-1}$. For $m$ even, the free involution on $\mathbb{C}P^{m-1}$ induces a free involution on $P(1,m-1)$. We feel that some exotic free involutions on $P(1,m-1)$ would possibly realize the cases (4) and (5) of the main theorem.
\end{itemize}

We conclude with the following remarks.

\begin{remark}\label{remark3}
For the three dimensional lens space $L_p^{3}(q)$, where $p= 4k$ for some $k$, Kim \cite[Theorem 3.6]{Kim1} has shown that the orbit space of any sense-preserving free involution on $L_p^{3}(q)$ is the lens space $L_{2p}^{3}(q')$, where $ q'q \equiv \pm 1$ or $q' \equiv \pm q$ mod $p$. Here an involution is sense-preserving if the induced map on $H_1(L_p^{3}(q); \mathbb{Z})$ is the identity map. This is the case (2) of the main theorem.
\end{remark}

\begin{remark}\label{remark4}
If $T$ is a free involution on $L_p^{3}(q)$ where $p$ is an odd prime, then $\mathbb{Z}_p$ and the lift of $T$ to $\mathbb{S}^3$ generate a group $H$ of order $2p$ acting freely on $\mathbb{S}^3$. The involution $T$ is said to be orthogonal if the action of $H$ on $\mathbb{S}^3$ is orthogonal. Myers \cite{Myers} showed that every free involution on $L_p^{3}(q)$ is conjugate to an orthogonal free involution. It is well known that there are only two groups of order $2p$, namely the cyclic group $\mathbb{Z}_{2p}$ and the dihedral group $D_{2p}$. But by Milnor \cite{Milnor}, the dihedral group cannot act freely and orthogonally on  $\mathbb{S}^3$. Hence $H$ must be the cyclic group $\mathbb{Z}_{2p}$ acting freely and orthogonally on $\mathbb{S}^3$. Therefore the orbit space $L_p^{3}(q)/ \langle T \rangle = \mathbb{S}^3/H = L_{2p}^{3}(q)$. Since $p$ is odd, $L_p^{3}(q)\simeq_2 \mathbb{S}^3$ and $L_p^{3}(q)/ \langle T \rangle \simeq_2 \mathbb{R}P^3$,  which is the case (1) of the main theorem.
\end{remark}
\bigskip

\section{An application to $\mathbb{Z}_2$-equivariant maps}\label{section7}
Let $X$ be a paracompact Hausdorff space with a fixed free involution and let $\mathbb{S}^n$ be the unit $n$-sphere equipped with the antipodal involution. Conner and Floyd \cite{Conner} asked the following question.

\vspace*{2mm}
\noindent Question: For which integer $n$, is there a $\mathbb{Z}_2$-equivariant map from $\mathbb{S}^n$ to $X$, but no such map from $\mathbb{S}^{n+1}$ to $X$?
\vspace*{2mm}

In view of the Borsuk-Ulam theorem, the answer to the question for $X = \mathbb{S}^n$ is $n$. Motivated by the classical results of Lyusternik- Shnirel'man \cite{Lyusternik}, Borsuk-Ulam \cite{Borsuk}, Yang \cite{Yang1, Yang2, Yang3} and Bourgin \cite{Bourgin}, Conner and Floyd defined the index of the involution on $X$ as $$\textrm{ind}(X) = max \{ n ~|~ \textrm{there exist a}~ \mathbb{Z}_2 \textrm{-equivariant map}~ \mathbb{S}^n \to X \}.$$

It is natural to consider the purely cohomological criteria to study the above question. The best known and most easily managed cohomology class are the characteristic classes with coefficients in $\mathbb{Z}_2$. Let $w \in H^1(X/G; \mathbb{Z}_2 )$ be the Stiefel-Whitney class of the principal $G$-bundle $X \to X/G$. Generalizing the Yang's index \cite{Yang2}, Conner and Floyd defined $$\textrm{{co-ind}}_{\mathbb{Z}_2}(X)~ = ~\textrm{largest integer}~n~ \textrm{such that}~ w^n \neq 0.$$ Since $\textrm{{co-ind}}_{ \mathbb{Z}_2}(\mathbb{S}^n)$ = $n$, by \cite[(4.5)]{Conner}, we have $$\textrm{ind}(X)\leq \textrm{{co-ind}}_{\mathbb{Z}_2}(X).$$

Also, since $X$ is paracompact Hausdorff, we can take a classifying map $$c : X/G \to B_G$$ for the principal $G$-bundle $X \to X/G$. If $k: X/G \to X_G$ is a homotopy equivalence, then $\rho k : X/G \to B_G$ also classifies the principal $G$-bundle $X \to X/G$ and hence it is homotopic to $c$. Therefore it suffices to consider the map $$\rho^*: H^1(B_G; \mathbb{Z}_2 ) \to H^1(X_G; \mathbb{Z}_2).$$ The image of the Stiefel-Whitney class of the universal principal $G$-bundle $G \hookrightarrow E_G \longrightarrow B_G$ is the Stiefel-Whitney class of the principal $G$-bundle $X \to X/G$.

Let $X \simeq_2 L_p^{2m-1}(q)$ be a finitistic space with a free involution. The Smith-Gysin sequence associated to the principal $G$-bundle $X \to X/G$ shows that the Stiefel-Whitney class is non-zero.

In case (1), $x \in H^1(X/G; \mathbb{Z}_2)$ is the Stiefel-Whitney class with $x^{2m}=0$. This gives $\textrm{{co-ind}}_{\mathbb{Z}_2}(X)=2m-1$ and hence $\textrm{ind}(X)\leq 2m-1$. Therefore, in this case, there is no $\mathbb{Z}_2$-equivariant map from $\mathbb{S}^n \to X$ for $n \geq 2m$.

Taking $X= \mathbb{S}^k$ with the antipodal involution, by proposition \ref{proposition5.1}, we obtain the classical Borsuk-Ulam theorem, which states that: There is no map from $\mathbb{S}^n \to \mathbb{S}^k$ equivariant with respect to the antipodal involutions when $n \geq k+1$.

In case (2), $x \in H^1(X/G; \mathbb{Z}_2)$ is the Stiefel-Whitney class with $x^2=0$. This gives $\textrm{{co-ind}}_{\mathbb{Z}_2}(X)=1$ and $\textrm{ind}(X)\leq 1$. Hence, there is no  $\mathbb{Z}_2$-equivariant map from $\mathbb{S}^n \to X$ for $n \geq 2$.

In case (3), $x \in H^1(X/G; \mathbb{Z}_2)$ is the Stiefel-Whitney class with $x^3=0$. This gives $\textrm{{co-ind}}_{\mathbb{Z}_2}(X)=2$ and $\textrm{ind}(X)\leq 2$. Hence, there is no  $\mathbb{Z}_2$-equivariant map from $\mathbb{S}^n \to X$ for $n \geq 3$.

In case (4), $x \in H^1(X/G; \mathbb{Z}_2)$ is the Stiefel-Whitney class with $x^4=0$. This gives $\textrm{{co-ind}}_{\mathbb{Z}_2}(X)=3$ and hence $\textrm{ind}(X)\leq 3$. Hence, there is no $\mathbb{Z}_2$-equivariant map from $\mathbb{S}^n \to X$ for $n \geq 4$.

Finally, in case (5) of the main theorem, $x \in H^1(X/G; \mathbb{Z}_2)$ is the Stiefel-Whitney class with $x^5=0$. This gives $\textrm{{co-ind}}_{\mathbb{Z}_2}(X)=4$ and hence $\textrm{ind}(X)\leq 4$. In this case also, there is no $\mathbb{Z}_2$-equivariant map from $\mathbb{S}^n \to X$ for $n \geq 5$.

Combining the above discussion, we have proved the following Borsuk-Ulam type result.

\begin{theorem}
Let $m \geq 3$ and $X \simeq_2 L_p^{2m-1}(q)$ be a finitistic space with a free involution. Then there does not exist any $\mathbb{Z}_2$-equivariant map from $\mathbb{S}^n \to X$ for $n \geq 2m$.
\end{theorem}

We end by pointing out another possible application of the main theorem.

\begin{remark}
The Borsuk-Ulam theorem has been extended to the setting of fiber bundles by many authors. It is known as parametrization of the Borsuk-Ulam theorem and a general formulation of the problem for free involutions is as follow.

\vspace*{2mm}
\noindent Problem: Let $X \hookrightarrow E \to B$ be a fiber bundle and $E' \to B$ be a vector bundle such that $\mathbb{Z}_2$ acts fiber preservingly and freely on both $E$ and $E'-0$, where $0$ stands for the zero section of the bundle $E' \to B$. For a fiber preserving $\mathbb{Z}_2$-equivariant map $f:E \to E'$, the parametrized Borsuk-Ulam problem is to estimate the cohomological dimension of the zero set $\{x \in E~|~f(x)=0 \}$ of $f$.
\vspace*{2mm}

Dold \cite{Dold2} and Nakaoka \cite{Nakaoka} considered the problem when $X$ is a cohomology sphere. Using their method, Koikara and Mukerjee \cite{Koikara} proved a similar theorem for $X$ a product of spheres. Recently, the author \cite{msingh2} proved a parametrized Borsuk-Ulam theorem when $X$ is a cohomology projective space. The main ingredient in these results is the knowledge of the cohomology algebra of $X/\mathbb{Z}_2$. We expect that the results in this paper could be used for obtaining a parametrized Borsuk-Ulam type theorem for lens space bundles.
\end{remark}
\bigskip

\begin{acknowledgement}
The author would like to thank Professor Pedro Luiz Queiroz Pergher for encouraging to publish results of this paper. The author would also like to thank the referee for valuable comments which improved the paper. Finally, the author would like to thank the Department of Science and Technology of India for support via the INSPIRE Scheme IFA-11MA-01/2011 and the SERC Fast Track Scheme SR/FTP/MS-027/2010.
\end{acknowledgement}
\bigskip

\bibliographystyle{amsplain}

\begin{thebibliography}{10}
\bibitem {Ashraf} R. Ashraf, \textit{Singular cohomology rings of some orbit spaces defined by free involution on $\mathbb{C}P(2m + 1)$}, J. Algebra 324 (2010), 1212-1218.
\bibitem {Borsuk} K. Borsuk, \textit{Drei S\"{a}tze \"{u}ber die n-dimensionale euklidische Sph\"{a}re}, Fund. Math. 20 (1933), 177-190.
\bibitem{Bourgin} D. G. Bourgin, \textit{On some separation and mapping theorems}, Comment. Math. Helv. 29 (1955), 199-214.
\bibitem {Bredon2} G. E. Bredon, \textit{Introduction to Compact Transformation Groups}, Academic Press, New York, 1972.
\bibitem {Bredon3} G. E. Bredon, \textit{Sheaf Theory}, Second Edition, Springer-Verlag, New York, 1997.
\bibitem {Borel3} A. Borel et al., \textit{Seminar on Transformation Groups}, Annals of Math. Studies 46, Princeton University Press, 1960.
\bibitem {Conner} P. E. Conner and E. E. Floyd \textit{Fixed point free involutions and equivariant maps-I}, Bull. Amer. Math. Soc. 66 (1960), 416-441.
\bibitem {Dold1} A. Dold \textit{Erzeugende der Thomschen algebra $\mathfrak{N}_*$}, Math. Z. 65 (1956), 25-35.
\bibitem {Dold2} A. Dold, \textit{Parametrized Borsuk-Ulam theorems}, Comment. Math. Helv. 63 (1988), 275-285.
\bibitem {Dotzel} R. M. Dotzel, T. B. Singh and S. P. Tripathi \textit{The cohomology rings of the orbit spaces of free transformation groups on the product of two spheres}, Proc. Amer. Math. Soc. 129 (2000), 921-930.
\bibitem {Hatcher} A. Hatcher, \textit{Algebraic Topology}, Cambridge University Press, 2002.
\bibitem {Kim1} P. K. Kim, \textit{Periodic homeomorphisms of the 3-sphere and related spaces}, Michigan Math. J. 21 (1974), 1-6.
\bibitem {Koikara} B. S. Koikara and H. K. Mukerjee, \textit{A Borsuk-Ulam type theorem for a product of spheres}, Topology Appl. 63 (1995), 39-52.
\bibitem {Livesay} G. R. Livesay, \textit{Fixed point free involutions on the $3$-sphere},  Ann. of Math. 72 (1960), 603-611.
\bibitem {Lyusternik} L. Lyusternik and L. Shnirel'man, \textit{Topological methods in variational problems and their application to the differential geometry of surfaces}, Uspehi Matem. Nauk. 2 (1947), 166-217.
\bibitem {Mccleary} J. McCleary, \textit{A User's Guide to Spectral Sequences}, Cambridge Studies in Advanced Mathematics 58, Second Edition, Cambridge University Press, Cambridge, 2001.
\bibitem {Milnor} J. Milnor, \textit{Groups which act on $\mathbb{S}^n$ without fixed points}, Amer. J. Math. 79 (1957), 623-630.
\bibitem {Myers} R. Myers, \textit{Free involutions on lens spaces}, Topology 20 (1981), 313-318.
\bibitem {Nakaoka} M. Nakaoka, \textit{Parametrized Borsuk-Ulam theorems and characteristic polynomials}, Topological fixed point theory and applications (Tianjin, 1988), 155-170, Lecture Notes in Math. 1411, Springer, Berlin, 1989.
\bibitem {Rice} P. M. Rice, \textit{Free actions of $\mathbb{Z}_4$ on $\mathbb{S}^3$}, Duke Math. J. 36 (1969), 749-751.
\bibitem {Ritter} G. X. Ritter, \textit{Free $\mathbb{Z}_8$ actions on $\mathbb{S}^3$}, Trans. Amer. Math. Soc. 181 (1973),195-212.
\bibitem {Ritter2} G. X. Ritter, \textit{Free actions of cyclic groups of order $2^n$ on $\mathbb{S}^1 \times \mathbb{S}^2$}, Proc. Amer. Math. Soc. 46 (1974), 137-140.
\bibitem {Rubinstein} J. H. Rubinstein, \textit{Free actions of some finite groups on $\mathbb{S}^3$ - I}, Math. Ann. 240 (1979), 165-175.
\bibitem{Swan} R. G. Swan, \textit{A new method in fixed point theory}, Comment. Math. Helv. 34 (1960), 1-16.
\bibitem{tbsingh} H. K. Singh and T. B. Singh, \textit{On the cohomology of orbit space of free $\mathbb{Z}_p$-actions on lens spaces},  Proc. Indian Acad. Sci. Math. Sci. 117 (2007), 287-292.
\bibitem{msingh1} M. Singh, \textit{Orbit spaces of free involutions on a product of two projective spaces}, Results. Math. 57 (2010), 53-67.
\bibitem{msingh2} M. Singh, \textit{Parametrized Borsuk-Ulam problem for projective space bundles}, Fund. Math. 211 (2011), 135-147.
\bibitem {Tao} Y. Tao, \textit{On fixed point free involutions on $\mathbb{S}^1 \times \mathbb{S}^2$}, Osaka J. Math. 14 (1962) 145-152.
\bibitem {Yang1} C. T. Yang, \textit{Continuous functions from spheres to Euclidean spaces}, Ann. of Math. 62 (1955), 284-292.
\bibitem {Yang2}  C. T. Yang, \textit{On theorems of Borsuk-Ulam, Kakutani-Yamabe-Yujobo and Dyson-I}, Ann. of Math. 60 (1954), 262-282.
\bibitem {Yang3}  C. T. Yang,\textit{On theorems of Borsuk-Ulam, Kakutani-Yamabe-Yujobo  and Dyson-II}, Ann. of Math. 62 (1955), 271-283.
\end{thebibliography}

\end{document}